\documentclass[11pt]{article}

\usepackage{graphicx, color}
\usepackage{amsmath,amsfonts,amssymb}
\def\forall{\hbox{for all}~}
\def\L{\mathbf{L}}

\def\ve{\varepsilon}

\def\D{{\cal D}}

\def\R{{\mathbb R}}

\def\vp{\varphi}

\def\tv{\hbox{Tot.Var.}}

\def\v{\vskip 1em}
\def\O{{\cal O}}
\def\C{{\cal C}}
\def\div{\hbox{div}}

\def\sign{{\rm sign}}
\def\bega{\begin{array}}
\def\enda{\end{array}}
\def\begi{\begin{itemize}}
\def\endi{\end{itemize}}

\def\ds{\displaystyle}
\def\bel{\begin{equation}\label}
\def\eeq{\end{equation}}
\def\sqr#1#2{\vbox{\hrule height .#2pt
\hbox{\vrule width .#2pt height #1pt \kern #1pt
\vrule width .#2pt}\hrule height .#2pt }}
\def\square{\sqr74}
\def\endproof{\hphantom{MM}\hfill\llap{$\square$}\goodbreak}

\newtheorem{theorem}{Theorem}[section]

\newtheorem{lemma}{Lemma}[section]

\newtheorem{remark}{Remark}[section]

\begin{document}
\title{\bf On Traffic Flow with Nonlocal Flux: A Relaxation  Representation}

\author{Alberto Bressan and Wen Shen\\ \, \\
Department of Mathematics, Penn State University.\\
University Park, PA~16802, USA.\\
\\
e-mails:~axb62@psu.edu,  
~wxs27@psu.edu}
\maketitle

\begin{abstract}
We consider a conservation law model of traffic flow, where the velocity of each car depends on a weighted average of the traffic density $\rho$ ahead. The averaging kernel is of exponential type: $w_\varepsilon(s)=\varepsilon^{-1} e^{-s/\varepsilon}$. By a transformation of coordinates, the problem can be reformulated as a  $2\times 2$ hyperbolic system with relaxation. Uniform BV bounds on the solution are thus obtained, independent of the scaling parameter $\varepsilon$. Letting $\varepsilon\to 0$, the limit yields a weak solution to the corresponding conservation law $\rho_t + ( \rho v(\rho))_x=0$. In the case where the velocity $v(\rho)= a-b\rho$ is affine,  using the Hardy-Littlewood rearrangement inequality  we prove that the limit is the unique entropy-admissible solution to the scalar conservation law.
\end{abstract}

\section{Introduction}
\label{sec:1}
\setcounter{equation}{0}

We consider a nonlocal PDE model for traffic flow, where
the traffic density $\rho=\rho(t,x)$ 
satisfies a scalar conservation law with nonlocal flux
\bel{1} \rho_t+(\rho v(q))_x~=~0.\eeq
Here 
$\rho\mapsto v(\rho)$ is a decreasing function, modeling the velocity of cars
depending on the traffic density, while the integral
\bel{2}q(x)~=~\int_0^{+\infty} w(s)\, \rho(x+s)\,ds\eeq
computes a weighted average of the car density.
On the function $v$ and the averaging kernel $w$, we shall always assume
\begi
\item[{\bf (A1)}] {\it The function $v:[0,\rho_{jam}]\mapsto \R_+$ is $\C^2$,
and satisfies}
\bel{a1}
v(\rho_{jam}) \,=\, 0,\qquad\qquad v'(\rho)~\leq~-\delta_*~<~0, \qquad\forall \rho\in [0, \rho_{jam}].
\eeq
%
\item[{\bf (A2)}]  {\it 
The  weight function $w\in \C^1(\R_+)$ satisfies
\bel{3}w'(s)\,\leq \,0, \qquad\int_0^{+\infty} w(s)\, ds~=~1.\eeq}
\endi

In \textbf{(A1)} 
one can think of  $\rho_{jam}$ as the maximum possible density of cars along the road, when all cars are packed bumper-to-bumper and nobody moves.
At a later stage, more specific choices for the functions $w$ and $v$ will be made.
In particular, we shall focus on the case  where $w(s) = e^{-s}$.

The conservation equation (\ref{1}) will be solved with initial data
\bel{id}
\rho(0,x)~=~\bar \rho(x)~\in~[0,\rho_{jam}]\,.\eeq

Given a weight  function $w$ satisfying (\ref{3}), we also consider the rescaled
weights
\bel{wep}w_\ve(s)\,\doteq\, \ve^{-1} w(s/\ve)\,.\eeq
As $\ve\to 0+$, the weight $w_\ve$ converges to a Dirac mass at the origin, 
and~the nonlocal equation \eqref{1} 
formally converges to the  scalar conservation law
\bel{Claw}
\rho_t + f(\rho)_x~=~0,\qquad\mbox{where} \qquad  f(\rho)\; \dot =\; \rho v(\rho). 
\eeq
The main purpose of this paper is to analyze the convergence of solutions of 
the nonlocal equation (\ref{1})  to those of (\ref{Claw}).

Conservation laws with nonlocal flux have attracted much interest in recent years,
because of their numerous applications and  the analytical challenges they pose. 
Applications of nonlocal models include  sedimentation~\cite{BBKT}, 
pedestrian flow and crowd dynamics~\cite{AggarwalGoatin2016, ColomboLecureuxMercier2011, ColomboGaravelloLecureuxMercier2011, ColomboGaravelloLecureuxMercier2012},
traffic flow~\cite{BG2016, CCS2019},
synchronization of oscillators~\cite{AmadoriHaPark2017}, 
slow erosion of granular matter~\cite{AmadoriShen2012},
materials with fading memory~\cite{ChenChristoforou2007}, 
some biological and industrial models~\cite{ColomboMarcelliniRossi2016},
and many others.
Due to the nonlocal  flux, the equation~\eqref{1}
behaves very differently from the classical conservation law \eqref{Claw}.
Its analysis faces additional difficulties and requires novel techniques. 

For a fixed weight function $w$, the
well posedness of the nonlocal conservation laws was proved in~\cite{BG2016}
with a Lax-Friedrich type numerical approximation, in \cite{KPS} by the method of characteristics,
and in~\cite{FKG2018} using a Godunov type scheme. 
Traveling waves for related nonlocal models  have been recently studied 
in~\cite{CS2019, RidderShen2018,  ShenDDDE2017, ShenTR, ShenKarim2017}.
See also the results for several space dimensions~\cite{AggarwalColomboGoatin2015},
and other related results in~\cite{CrippaLecureuxMercier2013, Zumbrun1999}.

Up to this date, however, the nonlocal to local limit for~\eqref{1} as $\ve\to 0+$ 
has remained a challenging question.
Namely, is it true that
the solutions of the Cauchy problem $\rho_\ve$ of (\ref{1})-(\ref{2}), 
with averaging kernels $w_\ve$ in (\ref{wep}),
as $\ve\to 0+$  
converge to the entropy admissible solutions of~\eqref{Claw}?
The question was already posed in~\cite{AmorimColomboTeixeira2015}. 
For a general weight function $w(\cdot)$, whose support covers an entire neighborhood of the origin, a negative answer  is provided by the counterexamples in~\cite{CCS2019}.
On the other hand, the results in~\cite{CCS2019} do not apply to the physically relevant models 
where the velocity $v$ is a monotone decreasing function and each driver only 
takes into account the
density of traffic ahead (not behind) the car.
Indeed, existence and uniqueness result for this more realistic model 
are given in~\cite{BG2016, CG18}.
Furthermore, various numerical simulations~\cite{AmorimColomboTeixeira2015, BG2016}
suggest that the behavior of $\rho_\ve$ should be stable in the limit $\ve\to 0+$. 
See also~\cite{CCS2019P} for the effect of numerical viscosity in the study of this limit.
In the case of monotone initial data, a convergence result was recently proved in \cite{KP}.

The main goal of the present paper is to 
study the  limit behavior of solutions to~\eqref{1}, for the averaging kernel 
$w_\ve(s)=\ve^{-1} \exp(-s/\ve)$, as $\ve \to 0$.  
In this setting, we first show that~\eqref{1} can be treated as a $2\times2$
system with relaxation, in a suitable coordinate system.
This formulation allows us to obtain a uniform bound on the total variation, 
independent of $\ve$.  As $\ve\to 0$, a standard compactness argument yields the convergence
$\rho_\ve\to \rho$ in $\L^1_{loc}$,  for a weak solution $\rho$ of (\ref{Claw}).
Finally, in the case of a Lighthill-Whitham speed~\cite{MR0072606, Whitham} of the form
$v(\rho)= a - b\rho$, we prove that the limit  solution $\rho$ coincides with  
the unique entropy weak solution of~\eqref{Claw}.

The remainder of the paper is organized as follows.
Section~\ref{s:2} contains a short proof of global existence, uniqueness, and continuous
dependence on the initial data, for solutions to (\ref{1})-(\ref{2}) with 
$v,w$, satisfying {\bf (A1)-(A2)}.
For Lipschitz continuous initial data,
solutions are constructed locally in time, as the fixed point of a contractive transformation.
By suitable a priori estimates,
we then show that these Lipschitz solutions can be extended globally in time. 
In turn, the semigroup of Lipschitz solutions can be continuously extended
(w.r.t.~the $\L^1$ distance) to a domain containing all initial data with bounded variation.

Starting with Section~\ref{s:3}, we restrict our attention to exponential kernels:
$w_\ve(s)= \ve^{-1} e^{-s/\ve}$.   In this case, the  conservation law with nonlocal
flux can be reformulated as a hyperbolic system with relaxation.   
In Section~\ref{s:4}, by a suitable transformation
of independent and dependent coordinates, we  establish a priori BV estimates
which are independent of the relaxation parameter $\ve$.  
We assume here that the initial density is uniformly positive.
By a standard compactness argument, in Section~\ref{s:5} we construct the limit
of a sequence of solutions with averaging kernels $w_\ve$, as $\ve\to 0$.
It is then an easy matter to show that any such limit 
provides a weak solution to the conservation law
(\ref{Claw}).  A much deeper issue is whether this limit coincides 
with the unique entropy-admissible solution.  In Section~\ref{s:6} 
we prove that this is indeed true, in the special case where the velocity function is affine:
$v(\rho) = a-b\rho$.  This allows a detailed analysis
of the convex entropy $\eta(\rho)= \rho^2$. 
Using the Hardy-Littlewood rearrangement inequality~\cite{HLP, LL}, 
we show that the entropy production is $\leq \O(1)\cdot\ve$.
Hence, in the limit as $\ve\to 0$, this entropy is dissipated.

We leave it as an open question to understand whether the same result is valid
for more general velocity functions $v(\cdot)$. Say, for $v(\rho)= a - b\rho^2$.
Moreover, all of our techniques heavily rely on the fact that the averaging kernel $w(\cdot)$
is exponential.  It would be of much interest to understand what happens 
for different kind  of kernels.

\section{Existence of solutions} 
\label{s:2}
\setcounter{equation}{0}
In this section we consider the Cauchy problem for (\ref{1})-(\ref{2}),
for a given initial datum
\bel{5}
\rho(0,x)~=~\bar \rho(x).\eeq
We consider the domain
\bel{dom}
\D~\doteq~\Big\{ \rho\in \L^\infty(\R)\,;\quad \tv\{\rho\}<\infty,\quad \rho(x) \in[0,\rho_{jam}]
\quad\forall x\in\R\Big\}.\eeq

\begin{theorem}\label{t:exist} Under the assumptions {\bf (A1)} and {\bf (A2)}, there exists a unique 
semigroup $S:[0,+\infty[\,\times\D \mapsto\D$, 
continuous in $\L^1_{loc}$, such that each trajectory
$t\mapsto S_t\bar \rho$ is a weak solution to the Cauchy problem 
(\ref{1})-(\ref{2}), (\ref{5}).
\end{theorem}

{\bf Proof.} 
We first construct a family of Lipschitz solutions, and show that 
they depend continuously on time and on the initial data, in the $\L^1$ distance.
By an approximation argument, we then construct solutions for 
general BV data $\bar \rho\in \D$.
\v
{\bf 1.}  Consider the domain of Lipschitz functions
\begin{eqnarray}
\D_L& \doteq& \Big\{ \rho\in \D\,;~~ \inf_x \rho(x)>0,~~ \sup_x \rho(x)
 < \rho_{jam}\,, 
\nonumber \\
&& \qquad \qquad
|\rho(x)-\rho(y)|\leq L|x-y|~~\forall 
x,y\in\R\Big\}.
\label{domL}
\end{eqnarray}
For every initial datum $\bar \rho\in \D_L$,
 we will construct  a solution 
$t\mapsto \rho(t,\cdot)\in\D_{2L}$ as the unique fixed 
point of a contractive transformation,
on a suitably small time interval $[0, t_0]$.

Given any function $t\mapsto \rho(t,\cdot)\in \D_{2L}$,
consider the corresponding integral averages 
\bel{qq}q(t,x)
~ = \int_0^\infty w(s) \rho(t,x+s)\, ds\,.
\eeq
We observe that
$$q_x(t,x)~=~
 \int_0^\infty w(s) \,\rho_x(t,x+s)\, ds.$$
 Hence
\bel{qlip}\|q_x(t,\cdot)\|_{\L^\infty}~\leq~ \|\rho_x(t,\cdot)\|_{\L^\infty} ~\leq~2L\,.
\eeq
Moreover, an integration by parts yields
\[
q_{xx}(t,x)= \int_0^\infty w(s) \,\rho_{xx}(t,x+s)\, ds
=-w(0) \rho_x(t,x) - \int_0^\infty w'(s)\, \rho_x(t,x+s)\, ds\,,
\]
therefore
\begin{eqnarray}
\|q_{xx}(t,\cdot)\|_{\L^\infty} &\leq& w(0) \|\rho_x(t,\cdot)\|_{\L^\infty}  + \|w'\|_{\L^1}\cdot 
\|\rho_x(t,\cdot)\|_{\L^\infty}  \nonumber \\
& = & 2 w(0) \|\rho_x(t,\cdot)\|_{\L^\infty} ~\leq~4L w(0).
\label{qxx}
\end{eqnarray}

Consider the transformation $\rho\mapsto u=\Gamma(\rho)$, where
$u$ is the solution to the linear Cauchy problem
\bel{LCP} u_t + (v(q)u)_x~=~0,\qquad \quad u(0,x)~=~\bar \rho(x), \qquad t\in [0, t_0],\eeq
with $q$ as in (\ref{qq}).
In the next two steps we shall prove:
\begi
\item[(i)] The values  $\Gamma(u)$ remain  uniformly bounded in the $W^{1,\infty}$ norm.
\item[(ii)] The map $\Gamma: \D_{2L}\mapsto \D_{2L}$ is  contractive  w.r.t.~the $\C^0$ norm.
\endi
By the contraction mapping theorem, a unique fixed point will thus exist, 
providing the solution to (\ref{LCP}) on the time interval $[0, t_0]$.
\v
{\bf 2.} To fix the ideas, assume 
\bel{brb}
0~<\delta_0~\leq~\bar \rho(x)~\leq ~\rho_{jam} -\delta_0\,,\eeq for some $\delta_0$.
From the equation
\bel{u3}
u_t+v(q) u_x~=~-v'(q) q_x\,,\qquad\qquad u(0,x)\,=\,\bar\rho,\eeq
integrating along characteristics and using (\ref{qlip}), we obtain
\bel{u4} 
\delta_0 - t\cdot \|v'\|_{\L^\infty} \, 2L~\leq~u(t,x)~\leq ~\rho_{jam} -\delta_0+t\cdot \|v'\|_{\L^\infty} \, 2L.
\eeq
Choosing $t_0< \delta_0\cdot (\|v'\|_{\L^\infty} \, 2L)^{-1}$, the solution $u$  will thus remain
strictly positive and smaller than $\rho_{jam}$, for all $t\in [0, t_0]$.
\v
{\bf 3.} Differentiating the conservation law in (\ref{LCP}) we obtain
\bel{uxt}u_{xt} + v(q)u_{xx}~=~-2 v'(q)q_x\,u_x - [v''(q) q_x^2  + v'(q)q_{xx}]u.\eeq
Let $Z(t)$ be the solution to the ODE
$$\dot Z~=~aZ + b,\qquad\qquad Z(0) ~=~L\,,$$
where
\[
 a ~\doteq~ 2\|v'\|_{\L^\infty}\cdot 2L\,,
 \qquad
 b ~\doteq~
 \Big[4L^2  \|v''\|_{\L^\infty} + 4L w(0) \|v'\|_{\L^\infty}\Big] \cdot \rho_{jam}\,.
\]
Since $$\|u_x(0,\cdot)\|_{\L^\infty}~=~\|\bar\rho_x\|_{\L^\infty}~\leq~L\,,$$
in view of (\ref{uxt}) and the bounds (\ref{qlip})--(\ref{qxx}),
a comparison argument yields
\bel{uxb} \|u_x(t, \cdot)\|_{\L^\infty}~\leq ~Z(t).\eeq
In particular, for $t\in [0, t_0]$ with $t_0$ sufficiently small, we have
\bel{ux3} \|u_x(t, \cdot)\|_{\L^\infty}~\leq ~2L\,.\eeq
\v
{\bf 4.} 
Using the identity
\[
q_x(t,x) ~= ~-w(0) \rho(t,x) - \int_0^\infty w'(s) \rho(t,x+s)\; ds
\]
and recalling that $w'(s)\leq 0$,
one obtains the bound
\bel{qx2}
\| q_x(t,\cdot)\|_{\L^\infty}~ \le ~2 w(0) \| \rho(t,\cdot)\|_{\L^\infty}.
\eeq

Next, consider two functions $t\mapsto \rho_1(t,\cdot)$, 
$t\mapsto \rho_2(t,\cdot)$, both taking values inside
$\D_{2L}$.
Then, for all $t\in[0,t_0]$, the corresponding weighted averages $q_1,q_2$
satisfy
\bel{qqt}
\|q_1(t,\cdot)-q_2(t,\cdot)\|_{W^{1,\infty}}~\leq~ 
(1+2 w(0))
\cdot \sup_{\tau\in [0, t_0]} \|\rho_1(\tau,\cdot)-\rho_2(\tau,\cdot)
\|_{\L^\infty}\,. \eeq
By choosing $t_0>0$ small enough,  we claim that the corresponding solutions $u_1, u_2$ of (\ref{LCP}) satisfy
\bel{contr}\|u_1(t,\cdot)-u_2(t,\cdot)\|_{\L^\infty}
~\leq~{1\over 2}\,\sup_{\tau\in [0,t]}\|\rho_1(\tau,\cdot)- \rho_2(\tau,\cdot)\|_{\L^\infty}
\qquad\forall t\in [0, t_0]\,.\eeq
Indeed, consider a point $(\tau,y)$.  Call $t\mapsto x_i(t)$, $i=1,2$, the corresponding characteristics.
These solve the equations 
\bel{ce}
\dot x_i ~=~v(q_i(t,x_i(t))),\qquad\qquad  x_i(\tau) = y.\eeq
Hence, moving backward in time, we have
\begin{eqnarray*}&&\hspace{-1.5cm} 
-{d\over dt} |x_1(t)-x_2(t)| \\
&\leq& \Big|v(q_1(t,x_1(t)))-v(q_1(t,x_2(t)))\Big|+ \Big|v(q_1(t,x_2(t)))- v(q_2(t,x_2(t)))\Big|
\\
&\leq& \|v'\|_{\L^\infty} \|q_{1,x}\|_{\L^\infty}\cdot |x_1(t)-x_2(t)| + \|v'\|_{\L^\infty} \|q_1-q_2\|_{\L^\infty}\,.
\end{eqnarray*}
By (\ref{qlip}), the quantity $\|q_{1,x}(t,\cdot)\|_{\L^\infty}$ remains uniformly bounded. 
The distance $Z(t) \doteq |x_1(t)-x_2(t)|$ between the two characteristics thus satisfies
a differential inequality of the form
$$-{d\over dt} Z(t)~\leq~a_* Z(t) +b_* \|q_1(t,\cdot)-q_2(t,\cdot)\|_{\L^\infty},\qquad  Z(\tau)=0,$$
for some constants $a_*,b_*$.   This implies
\bel{x12e} |x_1(t)-x_2(t)|~\leq~\int_t^\tau e^{(t-s) a_*}  \cdot b_* \|q_1(s,\cdot)-q_2(s,\cdot)\|_{\L^\infty}\, ds\,.
\eeq
The values $u_i(\tau, y)$, $i=1,2$, can now be obtained by integrating along characteristics. Indeed,
$${d\over dt} u_i(t, x_i(t))~=~v'(q_i(t, x_i(t)))\cdot q_{i,x}(t, x_i(t))\cdot u_i(t, x_i(t)),
\qquad u_i(0, x_i(0))~=~\bar \rho(x_i(0)).$$
Thanks to the a priori bounds (\ref{qxx}) on $\|q_{i, xx}(t,\cdot)\|_{\L^\infty}$, using (\ref{x12e})
for any $\epsilon>0$ we can choose $t_0>0$ such that 
$$|u_1(\tau,y)-u_2(\tau,y)|~\leq~\epsilon\cdot\sup_{t\in [0,\tau]} \|q_1(t,\cdot)-q_2(t,\cdot)\|_{\L^\infty}\,,$$
for all  $\tau\in [0, t_0]$ and $y\in\R$. In view of (\ref{qqt}), this implies (\ref{contr}).
\v
{\bf 5.} 
By the contraction mapping principle, there exists a unique function
$t\mapsto \rho(t,\cdot)$ such that $\rho(t,\cdot)= u(t,\cdot)$ for all $t\in [0, t_0]$.
This fixed point of the transformation $\Gamma$ provides the unique solution to the Cauchy 
problem (\ref{1})-(\ref{2}) with initial data  (\ref{5}).  
\v
{\bf 6.}  
In this step we show that this solution can be extended to all times $t>0$.
This requires (i)  a priori
upper and lower bounds of the form
\bel{uul}0~<~\delta_0~\leq~\rho(t,x)~\leq ~\rho_{jam}-\delta_0\,,\eeq
independent of time, and (ii) a priori estimates on the Lipschitz constant, which should remain uniformly bounded on bounded intervals of time.
\v
To establish an upper bound on the solution $\rho(t,\cdot)$, $t\in [0, t_0]$,
we analyze its behavior along a characteristic.
Fix $\epsilon>0$.  Consider any point $(\tau,\xi)$ such that 
$$\rho(\tau,\xi)~\geq~ \sup_{x\in\R}\rho(\tau,x)-\epsilon.$$
At
 the point $(\tau,\xi)$ one has
\begin{eqnarray}
\rho_t + v(q)\rho_x&=&- \rho v'(q) q_x~=~-\rho(\tau,\xi) v'(q(\tau,\xi)) 
\cdot{\partial\over \partial\xi}\left[ \int_\xi^{+\infty} \rho(\tau, y) \, w
(y-\xi)\, dy\right]
\nonumber\\
&=& -\rho(\tau,\xi) v'(q(\tau,\xi)) 
\cdot\left[
-\rho(\tau,\xi)\,w(0) - \int_\xi^{+\infty} \rho(\tau, y) \, w'
(y-\xi)\, dy
\right]
\nonumber\\
&=&-\rho(\tau,\xi) v'(q(\tau,\xi)) 
\cdot \int_\xi^{+\infty}\bigl[ \rho(\tau,\xi)- \rho(\tau, y) \bigr]\, w'
(y-\xi)\, dy
\nonumber\\
&\leq &\rho_{jam} \cdot \max_{0\leq q\leq \rho_{jam}}
 |v'(q)|\cdot w(0)\cdot \epsilon~\doteq~C_0\, \epsilon.
\label{8}
\end{eqnarray}
The above implies
$${d\over dt} \left(\sup_x \rho(t,x)\right)~\leq~C_0\epsilon ,$$
as long as $0< \rho(t,y)<\rho_{jam}$ for all $y\in\R$.

Since $\bar \rho$ satisfies (\ref{brb}) and $\epsilon>0$ is arbitrary, this establishes the upper bound in (\ref{uul}).
The lower bound is proved in an entirely similar way.
\v
Next, from the analysis in step {\bf 3} it follows
\bel{ul4}\|\rho_x(t,\cdot)\|_{\L^\infty} ~\leq~Z(t),\eeq
which immediately yields the a priori bound on the Lipschitz constant.
\v
By induction, we can thus construct a unique solution $\rho=\rho(t,x)$ 
on a sequence of time intervals $[0, t_0]$, $[t_0, t_1]$, $[t_1, t_2], ~\ldots$,
where the length of each interval $[t_k\, t_{k+1}]$ depends only on
(i) the constant $\delta_0$ in (\ref{uul}), and (ii) the Lipschitz constant of $\rho(t_k,\cdot)$.
Thanks to (\ref{ul4}), this Lipschitz constant remains $\leq Z(t_k)$.  This implies
$t_k\to +\infty$ as $k\to \infty$, hence the solution can be extended to all times $t>0$.

We remark that, by a further differentiation of the basic equation (\ref{1}), one can prove that, if $\bar \rho\in C^k$, then 
every derivatives up to order $k$ remains uniformly bounded on bounded intervals of time.
\v
{\bf 7.} To complete the proof, it remains to 
show that the semigroup of solutions can be extended by continuity to all initial data
 $\bar \rho\in \D$.

Toward this goal, we first  prove that the total variation of the solution $\rho(t,\cdot)$ remains
uniformly bounded on bounded time intervals.
Indeed, from
$$\rho_{xt} + (v(q) \rho_x)_x~=~- (v'(q)q_x\rho)_x\,,$$
it follows
\begin{eqnarray}
&& \hspace{-1.5cm} 
{d\over dt} \|\rho_x\|_{\L^1}~ \leq~\|(v'(q)q_x\rho)_x\|_{\L^1}
\nonumber \\
&\leq& \|v'\|_{\L^\infty} \|q_x\|_{\L^\infty} \|\rho_x\|_{\L^1} +
\|v'\|_{\L^\infty} \|q_{xx}\|_{\L^1}  \|\rho\|_{\L^\infty} +\|v''\|_{\L^\infty} \|q_x\|_{\L^\infty}
\|q_x\|_{\L^1}\|\rho\|_{\L^\infty}
\nonumber\\ 
& \leq& C \|\rho_x\|_{\L^1}\,.
\label{tvloc}
\end{eqnarray}
Above we used the estimates
\bel{qxqxx2}
\|q_x\|_{\L^1} \le \|\rho_x\|_{\L^1}, \qquad 
\|q_{xx}\|_{\L^1} \le 2 w(0) \cdot \|\rho_x\|_{\L^1}.
\eeq
Note that in \eqref{tvloc} the constant $C$ depends on the velocity function $v:[0, \rho_{jam}]\mapsto\R_+$ and the averaging kernel $w$,
but it does not depend on the Lipschitz constant $\|\rho_x\|_{\L^\infty}$ of the solution.
According to (\ref{tvloc}), the total variation of the solution grows at most at an exponential rate.
In particular, it remains bounded on bounded intervals of time. 
\v
{\bf 8.}  Thanks to  the a priori bounds (\ref{tvloc}) on the total variation and  (\ref{uxb}) on the 
Lipschitz constant, the solution can be extended to an arbitrarily large time interval $[0,T]$.
This already defines a family of trajectories $t\mapsto S_t\bar \rho$ defined for every $L>0$, every
$\bar \rho\in \D_L$, and $t\geq 0$.   

In order to extend the semigroup $S$ by continuity 
to the entire domain $\D$, we need to prove that  for every $t>0$ the map $\bar \rho\mapsto S_t\bar \rho$
is Lipschitz continuous w.r.t.~the $\L^1$ distance.

Indeed, consider a family of smooth solutions, say $\rho^\theta(t,\cdot)$, $\theta>0$.
Define the first order perturbations 
$$\zeta^\theta(t,\cdot)~=~\lim_{h\to 0} {\rho^{\theta+h}(t,\cdot) - \rho^\theta(t,\cdot)\over h}\,,
\qquad \quad Q^\theta(t,\cdot)~=~\lim_{h\to 0} {q^{\theta+h}(t,\cdot) - q^\theta(t,\cdot)\over h}\,.$$
Notice that 
$$Q^\theta(t,x)~=~\int_0^{+\infty} w(s) \,\zeta^\theta(t,x+s)\, ds.$$
Then  $\zeta^\theta$  satisfies the linearized equation
\bel{zeq}
\zeta_t+(v(q) \zeta)_x + \bigl(v'(q) Q \rho\bigr)_x ~=~0,\eeq
where for simplicity we dropped the upper indices.
Using the estimates
\bel{QQx}
\|Q(t,\cdot)\|_{\L^1} \le \|\zeta(t,\cdot)\|_{\L^1} ,\qquad
\|Q_x(t,\cdot)\|_{\L^1} \le 2 w(0) \cdot \|\zeta(t,\cdot)\|_{\L^1} \;,
\eeq
\bel{QQx2}
\|q_x(t,\cdot)\|_{\L^\infty} \le 2 w(0) \cdot \rho_{jam},\qquad 
\| Q(t,\cdot)\|_{\L^\infty} \le w(0) \cdot \| \zeta (t,\cdot)\|_{\L^1}\;,
\eeq
we compute
\begin{eqnarray} && \hspace{-1.5cm}
{d\over dt}\|\zeta(t,\cdot)\|_{\L^1}~\leq~\|(v'(q) Q\rho)_x\|_{\L^1}
\nonumber \\
& \leq& \|v''\|_{\L^\infty} \|q_x\|_{\L^\infty}\|Q\|_{\L^1}\|\rho\|_{\L^\infty}
+\|v'\|_{\L^\infty} \|Q_x\|_{\L^1}\|\rho\|_{\L^\infty} 
+  \|v'\|_{\L^\infty} \|Q\|_{\L^\infty} \|\rho_x\|_{\L^1}
\nonumber \\
& \leq& C(t)\cdot \|\zeta(t,\cdot)\|_{\L^1}\,.
\label{lipsem}
\end{eqnarray}
Here $C(t)$ depends on time because the total variation $\|\rho_x(t,\cdot)\|_{\L^1}$ may grow at an exponential rate. 
On the other hand, it is important to observe that $C(t)$ does not depend 
on the Lipschitz constant of the solutions. From (\ref{lipsem}) we deduce
\bel{z6}\|\zeta(t,\cdot)\|_{\L^1} ~\leq~\exp\left\{ \int_0^t C(\tau)\, d\tau\right\} \|\zeta(0,\cdot)\|_{\L^1}\,.\eeq
For any two Lipschitz 
solutions $\rho^0$, $\rho^1$ of (\ref{1})-(\ref{2}), we now construct a 1-parameter family of 
solutions 
$\rho^\theta(t,\cdot)$ with initial data
$$\rho^\theta(0, \cdot) ~=~
\theta \rho^1(0,\cdot) + (1-\theta)\rho^0(0,\cdot).$$
Using (\ref{z6}) one obtains
\begin{eqnarray} 
\|\rho^1(t,\cdot)-\rho^0(t,\cdot)\|_{\L^1}&\leq&
\int_0^1 \|\zeta^\theta(t,\cdot)\|_{\L^1}\, d\theta~\leq\int_0^1\exp\left\{ \int_0^t C(\tau)\, d\tau\right\} 
\cdot \|\zeta^\theta(0,\cdot)\|_{\L^1}\, d\theta
\nonumber \\
&\leq& 
\exp\left\{ \int_0^t C(\tau)\, d\tau\right\} \cdot \|\rho^1(0,\cdot)-\rho^0(0,\cdot)\|_{\L^1}\,.
\label{lip3}
\end{eqnarray}
This establishes Lipschitz continuity of the semigroup w.r.t.~the initial data.
Notice that this Lipschitz constant may well depend on time.
Since every initial datum $\bar\rho\in \D$ can be approximated in the $\L^1$ distance 
by a sequence of Lipschitz continuous
functions $\bar\rho_n\in \D_{L_n}$ (possibly with $L_n\to +\infty$), 
by continuity we obtain a unique semigroup defined on the entire domain $\D$.
\endproof

\begin{remark}\label{r:21}
 {\rm By the argument in step {\bf 6} of the above proof,  if the
initial condition satisfies
$$0~\leq ~\bar a~\leq ~\bar \rho(x)~\leq~ \bar b~\leq ~\rho_{jam}\qquad\forall x\in\R,$$
then the solution satisfies
$$\bar a~\leq ~\rho(t,x)~\leq ~\bar b\qquad\forall t\geq 0,~ x\in\R.$$
}
\end{remark}
\section{A hyperbolic system with relaxation}
\label{s:3}
\setcounter{equation}{0}

From now on, we focus on the  case where $w(s) = e^{-s}$,
so that the rescaled kernels are 
\[w_\ve(s)= \ve^{-1} e^{- s/\ve}.\]
This yields
\bel{41} {\partial\over \partial x }\left[ \int_x ^{+\infty} \rho(t, s) \, {1\over \ve}\, e^{-
(s- x )/\ve}\, ds\right]~
=~
-\, {1\over \ve}\,\rho(t, x ) +{1\over \ve}\, \int_ x ^{+\infty} \rho(t, s) \, {1\over \ve}\, e^{-
(s- x )/\ve}\, ds.\eeq
Therefore, the averaged density $q$ satisfies the ODE 
$$q_x~=~\ve^{-1} q - \ve^{-1} \rho\,.$$
The conservation law with nonlocal flux (\ref{1})-(\ref{2}) can thus be written as 
\bel{43}
\left\{\bega{cl} \rho_t + (\rho v(q))_x&=~0,\\[1mm]
q_x&\ds =~\ve^{-1} (q-\rho)\,.\enda\right.\eeq
To make further progress, 
we  choose a constant $K>v(0)$ and
consider new independent coordinates $(\tau,y)$ defined by
\bel{cc}
\tau \,= \, t- {x\over K},\qquad \qquad y\,=\,  x\,.
\eeq
For future use, we derive the relations between the partial derivative 
operators in these two sets of coordinates: 
\bel{vchan} \partial_\tau~=~\partial_t\,,\qquad  \partial_y~=~\partial_x +  K^{-1} \partial_t\,,\qquad \partial_x~=~\partial_y - K^{-1} \partial_\tau\,.\eeq
A direct computation yields  
\[
\rho_t = \rho_\tau ,\qquad  
(\rho v(q))_x = -K^{-1}  (\rho v(q))_\tau + (\rho v(q))_y,
\qquad q_x = -K^{-1}  q_\tau + q_y.\]
In these new coordinates, the equations (\ref{43})  take the form
\begin{equation}\label{w2b}
\left\{ \begin{array}{ccl}
\displaystyle
(K \rho -  \rho v(q))_\tau + (K \rho v(q))_y &=& 0, \\[2mm]
  q_\tau - K q_y &=&\displaystyle
 \frac{K}{\ve} (\rho-q).
\end{array}\right.
\end{equation}
One can easily verify that 
the above system of balance laws is strictly hyperbolic, with two 
distinct characteristic speeds
\bel{n1}
\lambda_1=-K, \qquad \lambda_2 = \frac{K v(q)}{K-v(q)} \,.
\eeq
We observe that
$\lambda_1 < 0 < \lambda_2$, provided that $K$ is sufficiently large such that $K>v(0)$. 
Moreover,  both characteristic families are linearly degenerate. 

In the zero relaxation limit, letting $\ve\to 0+$ one formally obtains 
$q\to \rho$.   Hence \eqref{w2b} formally converges to the scalar conservation law
\bel{n2}
(K \rho -  \rho v(\rho))_\tau + (K \rho v(\rho))_y ~=~ 0.
\eeq
Recalling  the function $f$ defined in \eqref{Claw}, 
one obtains
\bel{n5}
(K \rho -  f(\rho))_\tau + (K f(\rho))_y ~=~ 0.
\eeq
Note that \eqref{n5} is equivalent to the conservation law \eqref{Claw} 
in the original $(t,x)$ coordinates.

The characteristic speed for \eqref{n5} is
\bel{n3}
\lambda^* 
~=~\frac{K f'(\rho)}{K-f'(\rho)} ~=~ \frac{K^2}{K-f'(\rho)} -K.
\eeq
Since $K>v(0)\ge v(\rho) > f'(\rho)$, 
we clearly have $\lambda^* > -K=\lambda_1$.
Furthermore, since $f'(\rho) < v(\rho)$, we conclude that $\lambda^* < \lambda_2$.
The sub-characteristic condition 
\bel{sub}
\lambda_1\, < \,\lambda^* \,<\,\lambda_2
\eeq
is thus satisfied. This is a crucial condition for stability of the relaxation system, 
see \cite{MR0872145}.
For other related general references on zero relaxation limit, we refer to~\cite{ABWS2000, MR1213992}.

From \eqref{w2b} it follows
\begin{eqnarray*} 
(K-v(q))\rho_\tau + K v(q) \rho_y&=&\rho\bigl[ v(q)_\tau - K v(q)_y\bigr] \\
&=&\rho v'(q) (q_\tau - K q_y)~=~\rho v'(q)\cdot  \frac{K}{\ve} (\rho-q).
\end{eqnarray*}
We can thus write (\ref{w2b}) in  diagonal form: 
\begin{equation}\label{w11}
\left\{ \begin{array}{cl}
\displaystyle
\rho_\tau  + \frac{K v(q)}{K-v(q)} \rho_y & \displaystyle =~\frac{K}{\ve} \cdot (\rho-q) \cdot 
\frac{ \rho v'(q)}{K-v(q)},
\\[2.5mm]
\displaystyle
q_\tau - K q_y & \displaystyle =~ \frac{K}{\ve} \cdot  (\rho-q).
\end{array}\right.
\end{equation}
To further analyze (\ref{w11}), it is convenient to introduce the new dependent variables
\bel{uzdef}
u \,=\, \ln \rho, \qquad  z\,=\,\ln(K-v(q)),
\eeq
so that 
\bel{rhov}
\rho\,=\, e^u, \qquad v(q)\,=\, K-e^z  .
\eeq
Using these new variables, (\ref{w11}) becomes
\bel{w5}
\left\{
\begin{array}{cl}
\displaystyle 
u_\tau + K(Ke^{-z}-1) u_y &\displaystyle =~~ \frac{K}{\ve} \Lambda(u,z),\\[2.5mm]
\displaystyle z_\tau -K z_y &\displaystyle  =\, - \frac{K}{\ve} \Lambda(u,z),
\end{array}
\right.
\eeq
where the source term $\Lambda$ is given by
\bel{Ldef}
\Lambda(u,z) 
~=~  (\rho(u)- q(z)) \frac{v'(q(z))}{K-v(q(z))}\,.
\eeq
Introducing the monotone function 
\bel{g}
 g(u) ~\dot=~ \ln (K- v(e^u)), \qquad \quad
\mbox{where}\quad 
 g'(u) ~= ~\frac{-v'(e^u) e^u}{K-v(e^u)}~>~0,
\eeq
one checks that
\bel{g2}
\Lambda(u,g(u))\, =\,0 \qquad \forall u .
\eeq

Letting $\ve\to 0$, we expect that $z \to g(u)$ hence the system \eqref{w5} 
formally converges to the scalar conservation law
\bel{g3}
(u + g(u))_\tau + K(Ke^{-g(u)} -1) u_y - K g(u)_y =0.
\eeq
Using the identities
\[
 u+g(u)\, =\, \ln(e^u(K-v(e^u))),\quad 
 e^{-g(u)}\,=\,\frac{1}{K-v(e^u)},\quad
 Ke^{-g(u)} -1\, =\, \frac{v(e^u)}{K-v(e^u)}\,,
\]
we get
\[
\frac{(e^u(K-v(e^u)))_\tau}{e^u(K-v(e^u))} + \frac{K(e^u v(e^u))_y}{e^u(K-v(e^u))} ~=~0.
\]
Writing $\rho=e^u$,  we obtain once again  the conservation law \eqref{n2}.

\v

\section{A priori BV bounds}
\label{s:4}
\setcounter{equation}{0}

In order to prove  a rigorous convergence result, we need an
a priori BV bound on the solution to the system (\ref{w5}), 
independent of the relaxation parameter $\ve$.
We always assume that the velocity $v$ satisfies the assumptions {\bf (A1)}.

\begin{figure}[ht]
\centerline{\hbox{\includegraphics[width=11cm]{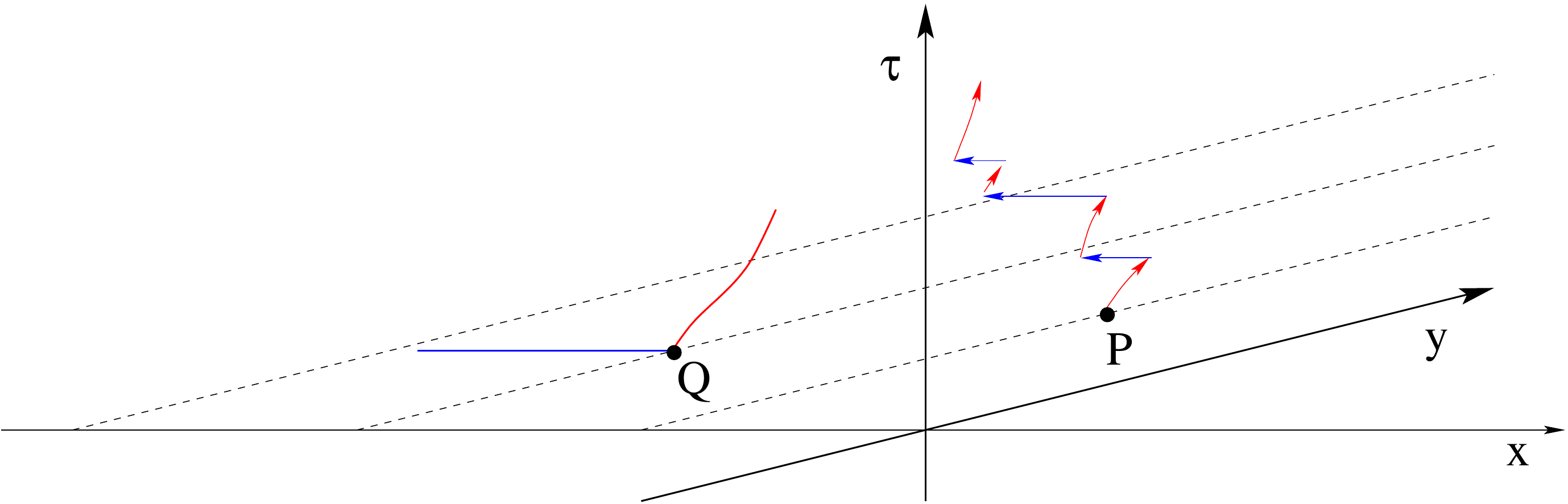}}}
\caption{\small The new system of coordinates $(\tau, y)$ defined at (\ref{cc}), 
is illustrated here together with the original coordinates $(t,x)$. 
The two characteristics through a point $Q$ have speeds 
$\lambda_1<0< 
\lambda_2$, as in (\ref{n1}).
With reference to the system (\ref{w6}), one can think of $z_y$ as the density
of backward-moving  particles, with speed $\lambda_1= -K$, while $u_y$ is the density
of forward-moving particles, with speed $\lambda_2>0$.  Backward particles are transformed
into forward particles at rate $K\Lambda_z/\ve$, while forward particles turn into 
backward ones with rate $-K\Lambda_u/\ve$.   The total number of particles does not increase;
actually, it decreases when positive and negative particles of the same type cancel out. }
\label{f:tf250}
\end{figure}

Differentiating (\ref{w5}) w.r.t.~$y$ one obtains
\bel{w6}
\left\{
\begin{array}{cl}
\displaystyle 
u_{y\tau} + [K(Ke^{-z}-1) u_y]_y &\displaystyle =~~\frac{K}{\ve}\,
[ \Lambda_u u_y + \Lambda_z z_y], \\[2.5mm]
\displaystyle z_{y\tau} -Kz_{yy} &\displaystyle  =~ -\frac{K}{\ve}\,
[ \Lambda_u u_y + \Lambda_z z_y] .
\end{array}
\right.
\eeq
A kinetic interpretation of the
above system  is shown in Figure~\ref{f:tf250}.

We observe that 
\begin{eqnarray*}
&& \hspace{-2cm} 
{d\over d\tau} \int |u_y(\tau,y)|\, dy + {d\over d\tau} \int |z_y(\tau,y)|\, dy \\
&=& \frac{K}{\ve}\,\int \Big\{ \sign(u_y) [ \Lambda_u u_y + \Lambda_z z_y] -\sign(z_y) [ \Lambda_u u_y + \Lambda_z z_y]\Big\}\, dy\\
&\leq & \frac{K}{\ve}\,\int  \Big\{\Lambda_u |u_y| + |\Lambda_z| \cdot |z_y| +   
|\Lambda_u| \cdot |u_y| - \Lambda_z |z_y|\Big\}\, dy.
\end{eqnarray*}
Therefore, if 
\bel{c4}\Lambda_u~\leq~0,\qquad\qquad \Lambda_z~\geq~0,\eeq
then the map
$$
\tau~\mapsto~\|u_y(\tau,\cdot)\|_{\L^1} + \|z_y(\tau,\cdot)\|_{\L^1}$$
will be non-increasing.
By (\ref{Ldef}), a direct computation yields
\[ \Lambda_u ~=~ e^u  \frac{v'(q(z))}{K-v(q(z))} ~< ~0.
\]

It remains to verify that $\Lambda_z \ge 0$.  
Since $\frac{\partial q}{\partial z} >0$, 
it suffices to show that $\Lambda_q \ge 0$.
We compute
\begin{eqnarray}
\Lambda_q &=& 
 (\rho-q) \frac{v''(q) (K-v(q)) +(v'(q))^2}{(K-v(q))^2}  - \frac{v'(q)}{K-v(q)}
\nonumber \\
&=& \frac{1}{K-v(q)} \left[ 
(\rho-q) \Big(v''(q) + \frac{(v'(q))^2}{K-v(q)} \Big) -v'(q)
\right].
\label{Lq}
\end{eqnarray}
Since  $v'(q)<0$, the above inequality will hold provided that
\bel{Lq2} |\rho-q|\cdot 
 \left( |v''(q)| + \frac{(v'(q))^2}{K-v(q)}\right)~ \le~ \left| v'(q) \right|\qquad \forall q.
\eeq
Notice that, by choosing $K$ sufficiently large, the factor ${(v')^2\over K-v}$ can be rendered as small as we like. Hence we can always achieve the inequality (\ref{Lq2})
provided that:
\begi
\item Either $|\rho-q|$ remains small.     This is certainly the case if
the oscillation of the initial datum is small.
\item Or else, $|v''|$ is small compared with $|v'|$.
\endi
As a consequence of  the  above analysis, we have: 
\begin{lemma}\label{l:3}
Let $(u,z)$ be a Lipschitz solution to the relaxation system (\ref{w6}).
Assume that $\rho(\tau,y)=e^{u(\tau,y)}\in [\rho_1,\rho_2]$ for all $(\tau,y)$, and moreover
\bel{bv3}\min_{q\in [\rho_1,\rho_2]} ~|v'(q)|~\geq~(\rho_2-\rho_1)
\cdot\left( \|v''\|_{\L^\infty} + {\|v'\|_{\L^\infty}^2\over K- \|v\|_{\L^\infty}}\right).
\eeq
Then the total variation function
\bel{TV1}
\tau~\mapsto~\|u_y(\tau,\cdot)\|_{\L^1(\R)} + \|z_y(\tau,\cdot)\|_{\L^1(\R)}\eeq
is non-increasing.
\end{lemma}

We observe that, 
in the case where $v$ is affine, say
\bel{ex1} v(\rho) ~=~a_o - b_o\rho\eeq
for some $a_o,b_o>0$, by (\ref{a1}) we can always choose $K$ large enough so that
\bel{bv3b}\rho_{jam}
\cdot  {\|v'\|_{\L^\infty}^2\over K- \|v\|_{\L^\infty}}
~\leq~\min_{0\leq q\leq \rho_{jam}} ~|v'(q)|\,.\eeq
Hence  (\ref{bv3}) is satisfied.
\v
Our main goal is to obtain uniform BV bounds for solutions to the nonlocal 
conservation law (\ref{1})-(\ref{2}).
This will be achieved by working in the $(\tau,y)$ coordinate system.

\begin{theorem}\label{t:41}  Consider the Cauchy problem for (\ref{1})-(\ref{2}),
with kernel $w(s)= \ve^{-1}e^{-s/\ve}$.  Assume that the velocity function $v$ satisfies
 \bel{bv6}\min_{\rho\in [0, \rho_{jam}]} ~|v'(q)|~>~\rho_{jam}
\cdot  \|v''\|_{\L^\infty([0, \rho_{jam}])}\,.\eeq
Moreover, assume that the initial density $\bar \rho$ has bounded variation and is
 uniformly positive.  Namely,
 \bel{rmin}
 0~<~\rho_{\min}~\leq~\bar \rho(x)~\leq \rho_{\max} \leq~\rho_{jam}\qquad\qquad\forall x\in\R.
 \eeq
 
Then the total variation remains uniformly bounded in time:
\bel{tv5}
\hbox{\rm Tot.Var.}\{\rho(t,\cdot)\}~\leq~ {\rho_{\max}\over\rho_{\min}}\cdot \hbox{\rm Tot.Var.}\{\bar\rho\}
\qquad\qquad \forall t\geq 0.\eeq
\end{theorem}
\v
{\bf Proof.} {\bf 1.}  Assume first that $\rho$ is Lipschitz continuous.
By (\ref{w6}) it follows
\bel{div2}
\div\begin{pmatrix} u_y\\[2mm] K(Ke^{-z}-1) u_y\end{pmatrix} +
\div\begin{pmatrix} z_y\\[2mm] -Kz_y\end{pmatrix}~=~0.\eeq
Thanks to (\ref{bv6}), we can choose a constant $K$ large enough so that (\ref{c4}) holds.
In this case we also have
\bel{div3}
\div\begin{pmatrix} |u_y|\\[2mm] K(Ke^{-z}-1) |u_y|\end{pmatrix} +
\div\begin{pmatrix} |z_y|\\[2mm] -K|z_y|\end{pmatrix}~\leq~0.
\eeq
In terms of the original $(t,x)$ coordinates, by (\ref{vchan})
the inequality (\ref{div3}) takes the form
\begin{eqnarray}
&& \hspace{-1cm} 
\partial_t\left( \Big| u_x + {u_t\over K}\Big| + \Big| z_x + {z_t\over K}\Big|\right)
+ \left( {1\over K} \partial_t + \partial_x\right) \left( K(Ke^{-z}-1) \Big| u_x + {u_t\over K}\Big|
- K\Big| z_x + {z_t\over K}\Big|\right)
\nonumber \\[1mm]
&=&
\partial_t\left( Ke^{-z}\Big| u_x + {u_t\over K}\Big| \right)
+ \partial_x \left( K(Ke^{-z}-1) \Big| u_x + {u_t\over K}\Big|
- K\Big| z_x + {z_t\over K}\Big|\right)~\leq~0.
\label{div4}
\end{eqnarray}
\v
{\bf 2.}
Integrating (\ref{div4}) over any time interval $[0,T]$, we obtain
\begin{eqnarray}
&& \hspace{-1.5cm} 
\int \frac{K}{K-v( q(T,x))}\Big| u_x(T,x) + {u_t(T,x)\over K}\Big| \, dx
\nonumber \\
&\leq&
\int  \frac{K}{K-v( q(0,x))} \Big| u_x(0,x) + {u_t(0,x)\over K}\Big| \, dx\,.
\label{v1}
\end{eqnarray}
Since we are choosing $K>v(0)\geq v(q(t,x))$ for all $t,x$, the above 
denominators remain uniformly positive and bounded.
This implies
\bel{v2}\int \Big| u_x(T,x) + {u_t(T,x)\over K}\Big| \, dx~\leq~C_K\cdot 
\int \Big| u_x(0,x) + {u_t(0,x)\over K}\Big| \, dx\,,\eeq
with $C_K\doteq {K\over K-v(0)}$.

{\color{black}
Repeating the same argument, with $K$ replaced by $\gamma K$ where $\gamma>1$, we obtain
\bel{v3}\int \Big| u_x(T,x) + {u_t(T,x)\over \gamma K}\Big| \, dx~\leq~C_{\gamma K}\cdot 
\int \Big| u_x(0,x) + {u_t(0,x)\over \gamma K}\Big| \, dx\,,\eeq
where the constant is now $C_{\gamma K} = {\gamma K\over \gamma K-v(0)}$.
\v

{\bf 3.}
Next, we observe that, for any two numbers $\alpha,\beta$ and any number $\gamma>1$ 
one has
$$\alpha~=~\frac{\gamma}{\gamma-1} \left(\alpha+{\beta\over \gamma}\right) - \frac{1}{\gamma-1} (\alpha+\beta),$$
so
$$
 |\alpha|~\leq~\frac{\gamma}{\gamma-1}  \left|\alpha+{\beta\over \gamma}\right| + \frac{1}{\gamma-1} |\alpha+\beta|.$$
Applying the above inequality with $\alpha= u_x$, $\beta= K^{-1}u_t$, from (\ref{v2})-(\ref{v3}) 
one obtains
\begin{eqnarray}
\int \bigl| u_x(T,x)\bigr| \, dx
&\leq&
\frac{\gamma C_{\gamma K}}{\gamma-1} 
\int \Big| u_x(0,x) 
+ {u_t(0,x)\over \gamma K}\Big| \, dx
\nonumber \\
&& \quad +~\frac{C_K}{\gamma-1} 
\int \Big| u_x(0,x) + {u_t(0,x)\over K}\Big| \, dx\,.
\label{v4}
\end{eqnarray}
\v
{\bf 4.} 
By the assumption  (\ref{rmin}) and Remark~\ref{r:21} it follows that $$0~<~\rho_{\min}~\leq~ \rho(t,x)~\leq~\rho_{\max}
\qquad
\forall ~t\geq 0,~x\in\R.$$
By the change of variables (\ref{uzdef})-(\ref{rhov}), one has
\bel{aru}
|u_x|~=~{|\rho_x|\over\rho}~\leq~{|\rho_x|\over\rho_{\min}}\,,\qquad 
|u_t|~=~{|\rho_t|\over\rho}~\leq~{|\rho_t|\over\rho_{\min}}\,,\qquad 
|\rho_x|~\leq~\rho_{\max} |u_x|\,.
\eeq
Combining (\ref{aru}) with (\ref{v4}) we conclude
\begin{eqnarray}
&& \hspace{-1cm} 
\rho_{\max}^{-1}  \int \bigl| \rho_x(T,x)\bigr| \, dx
 ~\leq~
  \int \bigl| u_x(T,x)\bigr| \, dx 
 \nonumber\\
 &\le& 
 \frac{\gamma C_{\gamma K}}{\gamma-1} 
\int \left| u_x(0,x) + {u_t(0,x)\over \gamma K}\right| \, dx~+~\frac{C_K}{\gamma-1}
\int \left| u_x(0,x) + {u_t(0,x)\over K}\right| dx
\nonumber \\
& \leq&\frac{\gamma C_{\gamma K}}{(\gamma-1)\rho_{\min}}
\int \left[| \rho_x(0,x)| + {|\rho_t(0,x)|\over \gamma K} \right] dx
\nonumber \\
&& \qquad +
{C_K\over(\gamma-1)\rho_{\min}}
\int \left[| \rho_x(0,x)| + {|\rho_t(0,x)|\over K} \right] dx.
\label{ar2}
\end{eqnarray}
We observe that
\begin{eqnarray*}
 \int |\rho_t(0,x)| \; dx &\leq& \int \left| \bar \rho_x(x)  v(q(0,x)) \right| + 
\left| \bar \rho(x) v'(q(0,x)) q_x(0,x) \right|\, dx\\
&\le& \|v\|_{\L^\infty} \cdot \|\bar\rho_x\|_{\L^1} + \rho_{\max} \cdot \|v'\|_{\L^\infty}\cdot  \|q_x(0,\cdot)\|_{\L^1}
~\le~ C_0  \|\bar\rho_x\|_{\L^1},
\end{eqnarray*}
where $C_0\,\doteq\,\|v\|_{\L^\infty} + \rho_{\max} \cdot\|v'\|_{\L^\infty}$ 
is a bounded constant.
Recalling the values of the constants $C_K, C_{\gamma K}$, from (\ref{ar2}) we obtain
\begin{eqnarray}
&& \hspace{-1.5cm} 
\int \bigl| \rho_x(T,x)\bigr| \, dx 
\nonumber \\
&\leq&{\rho_{\max}\over\rho_{\min}}
\cdot \frac{1}{\gamma-1}
\left( {\gamma^2 K\over \gamma K-v(0)}
\Big( 1+{ C_0\over \gamma K}
\Big) +  {K\over K-v(0)}\Big( 1+{ C_0\over K}\Big)\right)\cdot \|\bar\rho_x\|_{\L^1}\,.
\label{ar3}
\end{eqnarray}
Since the constant $K$ can be chosen arbitrarily large, letting $K\to +\infty$ in (\ref{ar3}) we obtain
\[ 
\| \rho_x(T,\cdot)\|_{\L^1}~\leq~ {\rho_{\max}\over\rho_{\min}} \cdot \frac{\gamma+1}{\gamma-1}
\cdot \|\bar\rho_x\|_{\L^1}\,.
\]
We note that as $K\to \infty$, (\ref{bv3}) reduces to (\ref{bv6}). 
Again, since $\gamma>1$ can be chosen arbitrarily large, letting $\gamma \to \infty$ we
obtain 
\bel{ar4} \| \rho_x(T,\cdot)\|_{\L^1}~\leq~ {\rho_{\max}\over\rho_{\min}} \cdot \|\bar\rho_x\|_{\L^1}\,.
\eeq}
For any Lipschitz solution,
this provides an a priori bound on the total variation, which does not
depend on time or on the relaxation parameter $\ve$.  
By an approximation argument we conclude that (\ref{tv5}) holds, 
for every uniformly positive  initial condition $\bar \rho$ with bounded variation.
\endproof

\section{Existence of a limit solution}
\label{s:5}
\setcounter{equation}{0}
Relying on the a priori bound on the total variation, proved in Theorem~\ref{t:41},
we now show the existence of a limit  $\rho= \lim_{\ve\to 0+}\rho_\ve$, 
which provides a weak solution to the conservation law (\ref{Claw}).

\begin{theorem}\label{t:51}
Let $\bar\rho: \R\mapsto [\rho_{\min}, \rho_{\max}]$ be a uniformly positive initial datum, with
bounded variation.
Call $\rho_\ve$ the corresponding solutions to (\ref{1})-(\ref{2}), with averaging kernel
$w_\ve(s)= \ve^{-1} e^{-s/\ve}$.  Then, by possibly extracting a subsequence $\ve_n\to 0$,
one obtains the convergence $\rho_{\ve_n}\to \rho$ in $\L^1_{loc}(\R_+\times\R)$. 
The limit function $\rho$ provides a weak solution to the conservation law (\ref{Claw}).
\end{theorem}

{\bf Proof.}  By Theorem~\ref{t:41}, all solutions $\rho_\ve(t,\cdot)$ have uniformly bounded total variation.   The same is thus true for the weighted averages
$q_\ve(t,\cdot)$, where
\bel{qep}q_\ve(t,x)~=~\int_0^{+\infty}\ve^{-1} e^{-s/\ve} \rho_\ve(t,x+s)\, ds\,.\eeq
By (\ref{1}), this implies that the map $t\mapsto \rho_\ve(t,\cdot)$ is uniformly 
Lipschitz continuous w.r.t.~the $\L^1$ distance.

By a compactness argument based on Helly's theorem (see for example Theorem~2.4 in~\cite{Bbook}), we can select a sequence $\ve_n
\downarrow 0$ such that
\begin{eqnarray}
\rho_{\ve_n}~\to~\rho&&\hbox{in} ~~\L^1_{loc}(\R_+\times \R),
\label{53} \\[1mm]
\rho_{\ve_n}(t,\cdot)~\to~\rho(t,\cdot)
&& \hbox{in} ~~\L^1_{loc}( \R),\qquad\hbox{for a.e.~} t\geq 0.
\label{54}
\end{eqnarray}
 By (\ref{qep}), it now follows
\begin{eqnarray*}
\bigl\|q_\ve(t,\cdot)- \rho_\ve(t,\cdot)\bigr\|_{\L^1}&=& \int\!\!\!\int_{x<y} \ve^{-1}e^{(x-y)/\ve} 
\bigl|\rho_\ve(t,y)-\rho_\ve(t,x)\bigr|\, dy\,dx \\
&\leq&
\int\!\!\!\int\!\!\!\int_{x<s<y} \ve^{-1}e^{(x-y)/\ve} 
\bigl|\rho_{\ve,x}(t,s)\bigr|\, ds\,dy\,dx\\
&=&
\int_{-\infty}^{+\infty}\left( \int_0^{+\infty} \int_0^{+\infty}
\,\ve^{-1} e^{-\sigma/\ve} e^{-\xi/\ve}\, d\xi\, d\sigma\right) \bigl|\rho_{\ve,x}(t,s)\bigr|\, ds\\[1mm]
&=& 
\ve\cdot \tv\{\rho_\ve(t,\cdot)\},
\end{eqnarray*} 
where the variables $\sigma=y-s$, $\xi = s-x$ were used.
Therefore, as $\ve_n\to 0$, we have the convergence
$q_{\ve_n}\to \rho$ in $\L^1_{loc}$.   By (\ref{1}), this implies that the limit function  $\rho=\rho(t,x)$ 
is a weak solution to the scalar conservation law
(\ref{Claw}).
\endproof

\section{Entropy admissibility of the limit solution}
\label{s:6}
\setcounter{equation}{0}

In the previous section we proved that, as $\ve\to 0$,
 any limit in $\L^1_{loc}$ of solutions 
$u_\ve$ to (\ref{1}), (\ref{id}) with $\bar \rho\in BV$ and $q_\ve$ given by (\ref{qep})
is a weak solution to the conservation law (\ref{Claw}).    A key question is whether
this limit is the unique entropy admissible solution.
The following analysis shows that this is indeed the case when the velocity function
is affine, namely
\bel{vaff}
v(\rho)~=~a-b\rho\,.\eeq

\begin{theorem}
\label{t:7} Let the velocity function $v$ be affine.  Consider any  uniformly 
positive initial
datum $\bar \rho\in BV$. Then as $\ve\to 0$, the corresponding
solutions $\rho_\ve$ to (\ref{1}), (\ref{qep}), (\ref{id}) converge to the unique entropy admissible solution of (\ref{Claw}).
\end{theorem}

{\bf Proof.} For simplicity, we consider the case where $v(\rho) = 1-\rho$.
The general case (\ref{vaff}) is entirely similar.
According to \cite{DLOW, Panov94}, to prove uniqueness
it suffices to prove that the limit solution dissipates one single strictly convex
entropy.
We thus consider the entropy and entropy flux pair
\bel{eep}\eta(\rho)\,=\,{\rho^2\over 2}\,,\qquad 
\qquad\psi(\rho)~=~{\rho^2\over 2}-{2\rho^3\over 3}\,.
\eeq
When $v(\rho) = 1-\rho$, the equation  \eqref{1} can be written as
\[ \rho_t + (\rho(1-\rho))_x ~=~ (\rho(1-\rho) - \rho(1-q))_x ~= ~(\rho (q-\rho))_x\,.
\]
Multiplying both sides by $\eta'(\rho)=\rho$, we obtain
\bel{ee5}
 \eta(\rho)_t +\psi(\rho)_x ~=~ \rho (\rho(q-\rho))_x ~ =~
  (\rho^2 (q-\rho))_x - (q-\rho) \rho \rho_x\,.
\eeq
Given a test function $\vp\in \C^1_c(\R)$,  $\vp\geq 0$,
we thus need to estimate the quantity 
\[ J ~=~ J_1 - J_2\,,\]
where
\begin{eqnarray} 
J_1 &\doteq& \int  (\rho^2 (q-\rho))_x \vp \; dx~ =~ - \int \rho^2 (q-\rho) \vp_x \; dx\,,
\label{J1}\\
J_2&\doteq&\int \bigl(q(x)-\rho(x)\bigr)\cdot \rho(x)\rho_x (x)\vp(x)\,dx
\nonumber \\
&=& \int \left(\int_x^{+\infty} {1\over\ve} e^{(x-y)/\ve}\left( \int_x^y\rho_x(s)\,ds\right) dy
\right) \rho(x)\,\rho_x(x)\, \vp(x)\,dx\,.
\label{t4}
\end{eqnarray}
Our ultimate goal is to show that 
$$J~\leq~\O(1)\cdot \ve.$$
Since we have 
\[
|J_1|~ \le ~\|\rho \|_{\L^\infty}^2 \cdot  \| \vp_x\|_{\L^\infty}  \cdot \int | q(x) -\rho(x)|\, dx ~= ~\O(1)\cdot \ve,
\]
it remains to show that
\begin{equation}\label{eJ2}
J_2 ~\geq~\O(1)\cdot \ve.
\end{equation}

A key tool to achieve this estimate is
\begin{lemma}{\bf (Hardy-Littlewood inequality).} {\it 
For any two functions $g_1,g_2\geq 0$ vanishing at infinity,
one has
\bel{HL}\int g_1(x) \, g_2(x) \, dx~\leq~\int g_1^*(x) g_2^*(x)\, dx,\eeq
where $g_1^*, g_2^*$ are the symmetric decreasing rearrangements of $g_1,g_2$, respectively.}
\end{lemma}
For a proof, see \cite{HLP} or \cite{LL}.

\v
Starting from (\ref{t4}) we compute
\begin{eqnarray*}
J_2&=&\int\!\!\int\!\!\int_{x<s<y}{1\over\ve} e^{(x-y)/\ve} \rho_x(s) \rho(x)\rho_x(x)\, \vp(x)\,  dy\, ds\, dx\\[1mm]
&=&\int\!\!\int_{x<s}  e^{(x-s)/\ve}\rho_x(s)\, \rho(x)\rho_x(x)\, \vp(x)\, dx \, ds
\\[1mm]
&=&\int  \left(\int_x^{+\infty} e^{-s/\ve}\rho_x(s)\, ds\right) e^{x/\ve}\rho(x)\rho_x(x)\, \vp(x)\,dx\\[1mm]
&=& - \int \rho^2(x) \rho_x(x)\, \vp(x) \,dx + {1\over \ve} \int\!\! \int_{x<s} e^{-s/\ve}\rho(s) \,
e^{x/\ve} \rho(x)\rho_x(x)\, \vp(x)\, dx\, ds \\[1mm]
&=& \int {\rho^3(x)\over 3} \vp_x(x)\, dx +{1\over  \ve} \int e^{-s/\ve} \rho(s)\
\left( \int_{-\infty}^s \Big({\rho^2(x)\over 2}\Big)_x\,e^{x/\ve}  \vp(x)\,
dx\right)  ds
\\[1mm]
&\dot=& A+B-C-D,
\end{eqnarray*}
where
\begin{eqnarray*} 
A &\dot=& \int {\rho^3(x)\over 3} \vp_x(x)\, dx\,, \\
B &\dot=&  {{1\over \ve} \int \rho(s)\,{\rho^2(s)\over 2}\,\vp(s)\, ds}\,,\\
C &\dot=& {1\over \ve^2} \int \int_{-\infty}^s e^{(x-s)/\ve} {\rho^2(x)\over 2}\, \rho(s)\, \vp(x)\, dx\, ds\,,\\
D&\dot=& 
{{1\over \ve} \int \int_{-\infty}^s e^{(x-s)/\ve} {\rho^2(x)\over 2}\, \rho(s)\, \vp_x(x)\, dx\, ds}\,.
\end{eqnarray*}

To achieve some cancellations, using a Taylor expansion of the term $C$ we obtain
\[
C \;\dot=\; C_1+C_2+C_3\,,
\]
where
\begin{eqnarray}
C_1& \dot=&{{1\over\ve^2} \int \int_{-\infty}^s e^{(x-s)/\ve} {\rho^2(x)\over 2}\, \rho(s)\, \vp(s)\, dx\, ds}\,,\nonumber\\
C_2& \dot=& {{1\over\ve^2}\int \int_{-\infty}^s e^{(x-s)/\ve} {\rho^2(x)\over 2}\, \rho(s)\, (x-s)\vp_x(x)\, dx\, ds}  \,,\nonumber
\\
C_3&\dot=& {
 {1\over\ve^2}\int \int_{-\infty}^s e^{(x-s)/\ve} {\rho^2(x)\over 2}\, \rho(s)\, {(x-s)^2\over 2}\vp_{xx}(\zeta)\, dx\, ds} \,.\label{t77}
\end{eqnarray}
In the integral for $C_3$, it is understood that for each $x,s$ 
one must choose a suitable $\zeta= \zeta(x,s)\in [x,s]$.

We now compare the  integrals $B$ and $C_1$.
 Without loss of generality one can assume 
$\vp = \phi^3$ for some $\phi\in \C^2_c$, $\phi\geq 0$.
For any $\sigma \geq 0$,
we now apply the  Hardy-Littlewood inequality with 
$$g_1(x)~=~\rho^2(x) \phi^2(x),\qquad g_2(x)~=~\rho(x+\sigma)\,\phi(x+\sigma),$$
and obtain
\bel{HL2} \int {\rho^2(x)\over 2}\,\rho(x)\,\vp(x)\, dx ~\geq~  \int  
 {\rho^2(x)\over 2}\, \phi^2(x)\cdot \rho(x+\sigma)\,  \phi(x+\sigma)\, dx\,.
\eeq
Indeed, the level sets of the two functions $\rho^2\phi^2$ and $\rho\phi$ are the same.
By (\ref{HL}), the integral on the right hand side of (\ref{HL2}) is maximum (and coincides with 
$\int g_1^* g_2^*\, dx$) when $\sigma=0$.

Performing the change of variable $s=x+\sigma$, 
a further integration w.r.t.~$s$ yields
\begin{eqnarray} 
B& =& \frac{1}{\ve}  \int {\rho^2(x)\over 2}\,\rho(x)\,\vp(x)\, dx 
~\ge~  \frac{1}{\ve}
\int   {\rho^2(x)\over 2}\, \phi^2(x)\cdot \rho(s)\,  \phi(s)\, dx \nonumber\\
& \geq& \frac{1}{\ve^2} \int \int_{-\infty}^s e^{(x-s)/\ve} {\rho^2(x)\over 2}\, \phi^2(x)\cdot \rho(s)\, \phi(s)\, dx\, ds ~\doteq~ B_1 - B_2\,,
\label{BHL}
\end{eqnarray}
where
\begin{eqnarray}
B_1 &\dot=&{1\over\ve^2} \int \int_{-\infty}^s e^{(x-s)/\ve} {\rho^2(x)\over 2}\, \rho(s)\, \phi^3(s)\, dx\, ds
~=~ C_1\,,
\label{vB1}\\
B_2 &\dot=&
 {1\over\ve^2} \int \int_{-\infty}^s e^{(x-s)/\ve} {\rho^2(x)\over 2}\, \phi(x)\rho(s)\, [\phi^2(s)- \phi^2(x)]\, dx\, ds.
 \nonumber
\end{eqnarray}

To compute the last integral for $B_2$ we use the Taylor expansion
$$\phi^2(s)- \phi^2(x)~=~2\phi(x)\phi_x(x) \cdot (s-x)+ [2\phi_x^2(\zeta) + 2\phi_{xx}(\zeta)]
\cdot {(s-x)^2\over 2}\,,$$
where $\zeta=\zeta(x,s)\in [x,s]$.
This yields
\begin{eqnarray*}
B_2&= & {1\over\ve^2} \int \int_{-\infty}^s e^{(x-s)/\ve} (s-x)\cdot {\rho^2(x)\over 2}\, \rho(s)\,\phi^2(x)\,2\phi_x(x) dx\, ds\\[1mm]
&& +~ {1\over\ve^2} \int \int_{-\infty}^s e^{(x-s)/\ve}{(s-x)^2\over 2} \cdot {\rho^2(x)\over 2}\,\rho(s)\, \phi(x)\,[2\phi_x^2(\zeta) + 2\phi_{xx}(\zeta)]
\, dx\, ds\\[1mm]
&=&B_{21}+B_{22}+B_{23}\,,
\end{eqnarray*}
where
\begin{eqnarray*}
B_{21} &\doteq&{1\over\ve^2} \int \int_{-\infty}^s e^{(x-s)/\ve} (s-x)\cdot \rho^3(x)\, \phi^2(x)\phi_x(x) dx\, ds \,, \\[1mm]
B_{22} &\doteq& {1\over\ve^2} \int \int_{-\infty}^s e^{(x-s)/\ve} (s-x)\cdot \rho^2(x)\, \left(\int_x^s\rho_x(\sigma)\, d\sigma\right)\,\phi^2(x)\phi_x(x) dx\, ds\,,
\\[1mm]
B_{23} &\doteq& {1\over\ve^2} \int \int_{-\infty}^s e^{(x-s)/\ve}{(s-x)^2\over 2} \cdot \rho^2(x)\,\rho(s)\, \phi(x)\,[\phi_x^2(\zeta) + \phi_{xx}(\zeta)]
\, dx\, ds \,.
\end{eqnarray*}

The  term  $B_{21}$ is computed by
\bel{B21}
B_{21} ~=~\int \rho^3(x)\, {\vp_x(x)\over 3}\, dx ~=~ A\,.
\eeq
Concerning $B_{22}$, 
using $\sigma, x$, and $\xi=s-x$ as variables of integration, we obtain
\begin{eqnarray}
 \left|B_{22}\right|&\leq& \|\rho\|_{\L^\infty}^2 \cdot {1\over 3}\|\vp_x\|_{\L^\infty}
\cdot  {1\over\ve^2}\int\!\!\int\!\!\int_{x<\sigma<s} e^{(x-s)/\ve}(s-x) |\rho_x(\sigma)|\, dx\,d\sigma\, ds \nonumber \\[1mm]
&=& \|\rho\|_{\L^\infty}^2 \cdot {1\over 3}\|\vp_x\|_{\L^\infty}
\cdot{1\over\ve^2} \int  \int_0^{+\infty} e^{-\xi/\ve} \xi \left(\int_{\sigma-\xi }^\sigma dx\right) d\xi
\,|\rho_x(\sigma)|\, d\sigma \nonumber \\[1mm]
&=& \|\rho\|_{\L^\infty}^2 \cdot {1\over 3}\|\vp_x\|_{\L^\infty}
\cdot  \int \left( \int_0^{+\infty} {e^{-\xi/\ve}\over\ve^2} \xi^2\, d\xi \right) |\rho_x(\sigma)|\, d\sigma
\nonumber \\
&=&\|\rho\|_{\L^\infty}^2 \cdot {1\over 3}\|\vp_x\|_{\L^\infty}
\cdot \|\rho_x\|_{\L^1}\cdot 2\ve\,.
\label{B22}
\end{eqnarray}

The term $B_{23}$ can be estimated by
\begin{eqnarray}
|B_{23}|
&  \leq& \|\rho\|^2_{\L^\infty}  \|\phi\|_{\L^\infty} \Big( \|\phi_x\|^2_{\L^\infty} + \|\phi_{xx}\|_{\L^\infty}\Big) 
\int |\rho(s)|\int_{-\infty}^s e^{(x-s)/\ve} {(x-s)^2\over 2\ve^2}\, dx\, ds
\nonumber \\
&=& \|\rho\|^2_{\L^\infty}   \|\phi\|_{\L^\infty}  \Big( \|\phi_x\|^2_{\L^\infty} +
 \|\phi_{xx}\|_{\L^\infty}\Big) \cdot \|\rho\|_{\L^1}
\int_0^{+\infty}e^{-\sigma/\ve} {\sigma^2\over 2\ve^2}\, d\sigma 
\nonumber \\
&=& \|\rho\|^2_{\L^\infty}   \|\phi\|_{\L^\infty}  \Big( \|\phi_x\|^2_{\L^\infty} +
 \|\phi_{xx}\|_{\L^\infty}\Big) \cdot \|\rho\|_{\L^1}\cdot  \ve\,.
\label{B23}
\end{eqnarray}

An entirely similar argument shows that
the integral defining $C_3$  at (\ref{t77}) also approaches zero as $\ve\to 0$.  Indeed,
\begin{eqnarray}
|C_3|&=&
 {1\over\ve^2}\left|\int \int_{-\infty}^s e^{(x-s)/\ve} {\rho^2(x)\over 2}\, \rho(s)\, {(x-s)^2\over 2}\vp_{xx}(\zeta)\, dx\, ds\right|
\nonumber \\
&\leq& \|\vp_{xx}\|_{\L^\infty} \cdot \|\rho\|^2_{\L^\infty} \cdot 
{1\over 2\ve^2} \int |\rho(s)|\int_{-\infty}^s e^{(x-s)/\ve} {(x-s)^2\over 2}\, dx\, ds
\nonumber\\
&\leq& \|\vp_{xx}\|_{\L^\infty} \cdot \|\rho^2\|_{\L^\infty} \cdot \|\rho\|_{\L^1}\cdot 
{1\over 2\ve^2} \int_0^{+\infty}e^{-\sigma/\ve} {\sigma^2\over 2}\, d\sigma 
\nonumber\\
& =&\|\vp_{xx}\|_{\L^\infty} \cdot \|\rho\|^2_{\L^\infty} \cdot \|\rho\|_{\L^1}\cdot {\ve\over 2}\,. 
\label{vC3}
\end{eqnarray}

Finally, we estimate the sum of the remaining two terms:
\begin{eqnarray*}
D+C_2
&=& {1\over \ve} \int \int_{-\infty}^s e^{(x-s)/\ve} {\rho^2(x)\over 2}\, \rho(s)\, \vp_x(x)\, dx\, ds
 \\
&&\ds\qquad -{1\over \ve^2} \int \int_{-\infty}^s e^{(x-s)/\ve} {\rho^2(x)\over 2}\, \rho(s)\, (s-x)\vp_x(x)\, dx\, ds
 \\
&=& \int {\rho^2(x)\over 2} \vp_x(x)\left(  \int_x^{+\infty} e^{(x-s)/\ve} 
\Big({1\over\ve}-{s-x\over \ve^2} 
\Big)
 \rho(s) \, ds \right)\, dx\,.
 \end{eqnarray*}
Using the identity
\[
 \int_x^{+\infty} e^{(x-s)/\ve} \Big({1\over\ve}-{s-x\over \ve^2} \Big) \, ds =0\,,
\]
we compute 
 \begin{eqnarray*}
D+C_2&=& \int {\rho^2(x)\over 2} \vp_x(x)\left(  \int_x^{+\infty} e^{(x-s)/\ve} 
\Big({1\over\ve}-{s-x\over \ve^2} 
\Big)
[ \rho(s)-\rho(x)] \, ds \right)\, dx
 \\
& =& \int {\rho^2(x)\over 2} \vp_x(x) \int_x^{+\infty} e^{(x-s)/\ve} \Big({1\over\ve}-{s-x\over \ve^2} 
\Big) \int_x^s\rho_x(\sigma)\, d\sigma \, ds \, dx
 \\
& =&  \int {\rho^2(x)\over 2} \vp_x(x) \int_x^{+\infty} \rho_x(\sigma)\, \int_\sigma^{+\infty}    e^{(x-s)/\ve} \Big({1\over\ve}-{s-x\over \ve^2} 
\Big)\, ds\, d\sigma\, dx 
 \\
& =&  \int {\rho^2(x)\over 2} \vp_x(x) \int_x^{+\infty} \rho_x(\sigma)\, 
 e^{(x-\sigma)/\ve} \, \frac{x-\sigma}{\ve}\,
 d\sigma\, dx 
  \\
&=& \int \rho_x(\sigma) \int_{-\infty}^\sigma{\rho^2(x)\over 2} \vp_x(x)
  e^{(x-\sigma)/\ve} \, \frac{x-\sigma}{\ve}\,dx\,
 d\sigma\,.
\end{eqnarray*}
As a consequence, we obtain the following estimate
\begin{eqnarray} 
|D+C_2|
&\leq& 
\|\rho_x\|_{\L^1} \cdot \Big\|{\rho^2\over 2}\Big\|_{\L^\infty} \cdot \|\vp_x\|_{\L^\infty}
\cdot  \int_{-\infty}^\sigma e^{(x-\sigma)/\ve} \, \frac{\sigma-x}{\ve}\,dx
\nonumber\\ 
&=& 
\|\rho_x\|_{\L^1} \cdot \Big\|{\rho^2\over 2}\Big\|_{\L^\infty} \cdot \|\vp_x\|_{\L^\infty}
\cdot\ve\,.
\label{ae1}
\end{eqnarray}

\v
Summarizing all the above estimates (\ref{t77})-(\ref{ae1}), we have
\begin{eqnarray} 
J_2
&=&A+B-C-D \nonumber \\
&\geq& A+B_1-(B_{21}+B_{22}+B_{23}) -(C_1+C_2+C_3)-D\label{t8} \\
&=&(A-B_{21})+(B_1-C_1)  - (D+C_2) - B_{22}- B_{23} - C_3 \label{t8-4}  \\
&=&\O(1)\cdot\ve. \nonumber
\end{eqnarray}
Indeed, on the line \eqref{t8-4} the first two terms are zero, while the remaining
four terms have size $\O(1)\cdot\ve$.
Letting $\ve\to 0$ we thus obtain the desired entropy inequality.

We remark that the inequality on the line (\ref{t8}), accounting for possible entropy dissipation, 
is due to the relation $B\geq B_1-B_2$ in (\ref{BHL}). This follows from
the Hardy-Littlewood rearrangement inequality.
\endproof


\begin{thebibliography}{99}


\bibitem{AggarwalColomboGoatin2015}
  \newblock A.~Aggarwal, R.~M.~Colombo,  and P.~Goatin.
  \newblock Nonlocal systems of conservation laws in several space dimensions.
  \newblock \emph{SIAM J.~Numer.~Anal.}, \textbf{53} {(2015)}, 963--983.
  

\bibitem{AggarwalGoatin2016}
  \newblock A.~Aggarwal  and P.~Goatin.
  \newblock Crowd dynamics through nonlocal conservation laws.
  \newblock \emph{Bull. Brazilian Math. Soc.}, \textbf{47} {(2016)}, 37--50.

  
\bibitem{AmadoriHaPark2017}
  \newblock D. Amadori, S.Y. Ha,   and J. Park.
  \newblock On the global well-posedness of BV weak solutions to the Kuramoto--Sakaguchi equation,
  \newblock \emph{J. Differential Equations}, \textbf{262} {(2017)}, 978--1022.

\bibitem{AmadoriShen2012}
  \newblock D. Amadori  and W. Shen.
  \newblock Front tracking approximations for slow erosion,
  \newblock \emph{Discr. Contin. Dyn. Syst.}, \textbf{32} {(2012)}, 1481--1502.


\bibitem{AmorimColomboTeixeira2015}
  \newblock P.~Amorim, R.~M.~Colombo,   and  A.~Teixeira.
  \newblock On the numerical integration of scalar nonlocal conservation laws.
  \newblock \emph{EASIM: M2MAN}, \textbf{49} {(2015)}, 19--37.
  
%
%
%
%

\bibitem{BBKT}
\newblock F.~Betancourt, R.~B\"{u}rger, K.~H.~Karlsen, and E.~M.~Tory. 
\newblock On nonlocal conservation laws modeling sedimentation. 
\newblock \emph{Nonlinearity}, \textbf{24}, 855--885, 2011.


\bibitem{BG2016}
  \newblock S.~Blandin and P.~Goatin.
  \newblock Well-posedness of a conservation law with nonlocal flux arising in traffic flow modeling.
 \newblock \textit{Numer.~Math.} \textbf{132} (2016), 217--241.


\bibitem{Bbook}
\newblock A.~Bressan,
\newblock \emph{Hyperbolic Systems of Conservation Laws.
The One Dimensional Cauchy Problem}, 
\newblock Oxford University Press, 2000.

\bibitem{ABWS2000}
\newblock A.~Bressan and W.~Shen,
\newblock BV estimates for multicomponent chromatography with relaxation, 
\newblock \textit{Discr. Cont. Dynam. Syst.}, The Millennium Issue, \textbf{6} (2000), pp. 21--38.

\bibitem{ChenChristoforou2007}
  \newblock G.-Q. Chen  and  C. Christoforou.
  \newblock Solutions for a nonlocal conservation law with fading memory,
  \newblock \textit{Proc. Amer. Math. Soc.} \textbf{135} (2007), 3905--3915. 

\bibitem{MR1213992} 
 \newblock G. Q. Chen and T. P. Liu, 
 \newblock Zero relaxation and dissipation limits for hyperbolic conservation laws, 
 \newblock \textit{Comm. Pure Appl. Math.} \textbf{46} (1993), 755--781.
 
 
 \bibitem{CG18}
 \newblock  F. A. Chiarello and P. Goatin. 
 \newblock Global entropy weak solutions for general non-local traffic flow models with anisotropic kernel. 
 \newblock \textit{ESAIM: Math. Mod. Numer. Anal.}, \textbf{52} (2018), 163--180.
 
%

\bibitem{CS2019}
\newblock J.~Chien and W.~Shen.
\newblock Stationary wave profiles for nonlocal particle models of traffic flow on rough roads.
\newblock \textit{Nonlin. Diff. Equat. Appl.}  (2019) \textbf{26}: 53.

\bibitem{CCS2019}
\newblock M.~Colombo, G.~Crippa,  and  L.V.~Spinolo. 
\newblock On the singular local limit for conservation laws with nonlocal fluxes.
\newblock \textit{Arch.~Rational~Mech.~Anal.} \textbf{233} (2019),  1131--1167.

\bibitem{CCS2018}
\newblock M.~Colombo, G.~Crippa,   and L.V.~Spinolo.
\newblock Blow-up of the total variation in the local limit of a nonlocal traffic model.
\newblock Preprint 2018, arxiv:1902.06970.

\bibitem{CCS2019P}
\newblock M.~Colombo, G.~Crippa, M.~Graff,   and L.V.~Spinolo.
\newblock On the role of numerical viscosity in the study of the local limit of nonlocal conservation laws.
\newblock Preprint 2019, arxiv:1902.07513.



\bibitem{ColomboLecureuxMercier2011}
  \newblock R.~M.~Colombo  and  M.~L\'{e}cureux-Mercier.
  \newblock Nonlocal crowd dynamics models for several populations.
  \newblock \emph{Acta Math.~Sci.},  \textbf{32} {(2012)}, 177--196.

\bibitem{ColomboGaravelloLecureuxMercier2011}
  \newblock R.~M.~Colombo, M.~Garavello,   and M.~L\'{e}cureux-Mercier.
  \newblock Nonlocal crowd dynamics.
  \newblock \emph{C. R. Acad. Sci.~Paris, Ser.~I}, \textbf{349} {(2011)}, 769--772.

\bibitem{ColomboGaravelloLecureuxMercier2012}
  \newblock R.~M.~Colombo, M.~Garavello,   and M.~L\'{e}cureux-Mercier.
  \newblock A class of nonlocal models for pedestrian traffic.
  \newblock \emph{Math.~Models Methods Appl.~Sci.}, \textbf{22} {(2012)}.

\bibitem{ColomboMarcelliniRossi2016}
  \newblock R.~M.~Colombo, F.~Marcellini,   and E.~Rossi.
  \newblock Biological and industrial models motivating nonlocal conservation laws: A review of analytic and numerical results.
  \newblock \emph{Netw.~Heterog.~Media}, \textbf{11} {(2016)}, 49--67.

%



\bibitem{CrippaLecureuxMercier2013}
  \newblock G.~Crippa  and  M.~L\'{e}cureux-Mercier.
  \newblock Existence and uniqueness of measure solutions for a system of continuity equations with nonlocal flow.
  \newblock \emph{Nonlin. Diff. Equat. Appl.}, \textbf{20} {(2013)}, 523--537.
  
%
%


\bibitem{DLOW}
\newblock C.~De Lellis, F.~Otto, and M.~Westdickenberg, 
\newblock  Minimal entropy conditions for Burgers equation. 
\newblock \emph{Quart. Appl. Math.} \textbf{62} (2004), 687--700.

%
%

%
%

%
%

%
%

\bibitem{FKG2018} 
\newblock    J.~Friedrich, O.~Kolb,   and S.~G\"{o}ttlich.
\newblock A Godunov type scheme for a class of LWR traffic flow models with non-local flux, \newblock \emph{Netw.~Heterog.~Media} \textbf{13} (2018), 531--547.

%
%


\bibitem{HLP}\newblock{G.~H.~Hardy, J.~E.~Littlewood, and G.~Polya.}
 \newblock \emph{Inequalities.} 
  \newblock Cambridge University Press,  1952. 


%
%

\bibitem{KP} \newblock A.~Keimer and L.~Pflug, 
\newblock{On approximation of local conservation laws by nonlocal conservation laws.}
\newblock\emph{J. Math. Anal. Appl.} {\bf 475} (2019), 1927--1955.



\bibitem{KPS} \newblock A.~Keimer, L.~Pflug, and M.~Spinola,
\newblock  Existence, uniqueness and regularity results on nonlocal balance laws. 
\newblock \emph{J. Differential Equations} {\bf  263} (2017), 4023--4069.

\bibitem{LL}
\newblock E.~Lieb and M.~Loss, {\it Analysis.} 
\newblock Second edition. American Mathematical Society, Providence, 2001.

    
    
\bibitem{MR0072606} 
     \newblock  M. J.  Lighthill  and  G. B.  Whitham.
     \newblock  On kinematic waves. II. A theory of traffic flow on long crowded
   roads,
     \newblock \emph{Proc. Roy. Soc. London. Ser. A.}, \textbf{229} (1955), 317--345.


\bibitem{MR0872145} 
\newblock T. P. Liu, 
\newblock Hyperbolic conservation laws with relaxation, 
\newblock \textit{Commun. Math. Phys.} \textbf{108} (1987), 153--175. 


\bibitem{Panov94}
   \newblock E.~Y.~Panov. 
      \newblock Uniqueness of the solution of the Cauchy problem for a first order quasilinear
equation with one admissible strictly convex entropy. 
\newblock (Russian)
 \emph{Mat. Zametki} \textbf{55} (1994), 116--129; 
 \newblock translation in \emph{Math. Notes} \textbf{55} (1994), 517--525.



 \bibitem{RidderShen2018}
  \newblock J.~Ridder  and W.~Shen. 
  \newblock Traveling waves for nonlocal models of traffic flow.
 \newblock   \textit{Discr. Contin.~Dyn.~Syst.} (2019), to appear.  
 
%


	
\bibitem{ShenDDDE2017}
\newblock W.~Shen.
\newblock Traveling wave profiles for a Follow-the-Leader model for traffic flow  with rough road condition, 
\newblock  \textit{Netw. Heterog. Media},  \textbf{13} (2018), 449--478.  


\bibitem{ShenTR}
\newblock W.~Shen.
\newblock Traveling waves for conservation laws with nonlocal flux for traffic flow on 
rough roads,
\newblock \textit{Netw. Heterog. Media} (2019), to appear.


\bibitem{ShenKarim2017}
\newblock {W.~Shen and K.~Shikh-Khalil}.
 \newblock  {Traveling waves for a microscopic model of traffic flow},
 \newblock  \textit{Discr.~Cont.~Dyn.~Syst. - A}, \textbf{38} (2018), 2571--2589.
 




\bibitem{Whitham}
\newblock B.~Whitham.
\newblock \textit{Linear and Nonlinear Waves}.
\newblock Wiley \& Sons, New York, 1974. 



\bibitem{Zumbrun1999}
 \newblock K.~Zumbrun.
 \newblock {On a nonlocal dispersive equation modeling particle suspensions.}
 \newblock \emph{Quart. Appl. Math.}, \textbf{57} {(1999)}, 573--600.





\end{thebibliography}
\end{document}